# Pseudodiohpantine geometric figure. In the wake of one task of Hugo Steinhaus

**Zurab Aghdgomelashvili**
Doctor of Mathematics

**Department of Mathematics, Georgian Technical University, 77, M. Koctava str. 0160, Tbilisi, Georgia**

E-mail: z.aghdgomelashvili@gtu.ge,

**Abstract**. In the work is considered one of up to now unsolved by Hugo Steinhaus task on having integer length square and located in its plane point that are on integer distances from its vertexes.

**Introduction**

Hugo Steinhaus has several interesting tasks. Let's consider one of them.

Let us consider the still unsolved task of Hugo Steinhaus.

**Task 1**. (Task □) "Is it possible to construct any square with sides of integer length and indicate in its plane such a point M, from that the distances to all four vertices are expressed in natural numbers?" In the same place, Hugo Steinhaus made the remark "The solution of this most complicated tasks I do not know?"

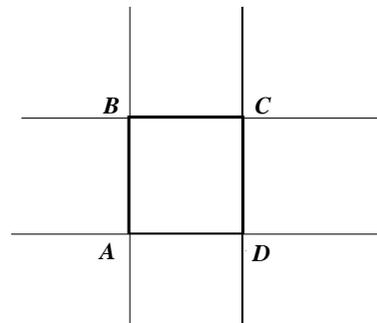

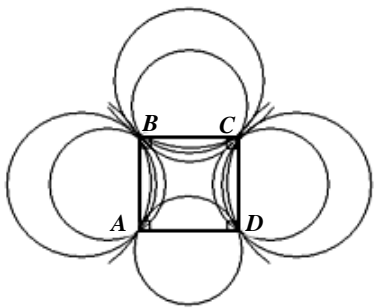

Hugo Steinhaus, as would be seen from one of this tasks (Task 1, does exist or not such a triangle, the length of each side of that is expressed by a natural number and any height of it is equal to the corresponding base?), by solving this task, he proved that if does exist such a square and if does exist such a point *M* in its plane, then the point *M* would not be located on the lines containing the sides of this square.

In this regard, we have obtained very good results. In this work it is shown that if does exist such a square and a point *M*, then it would not be located on such circumferences, the chords of that are represented by any side of this square, due that the circumference is divided into commensurate to π arcs.

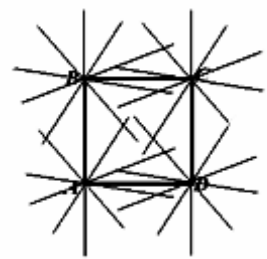

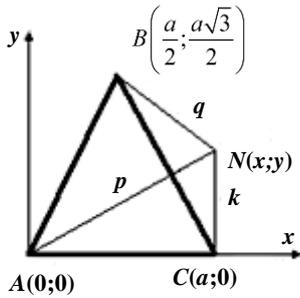

It is also shown that if such a point does exist, then it is impossible to be located on such a straight line passing through any vertex of this square that makes an angle commensurate with π with the straight lines containing the sides of the square.

We have set a similar Task for isosceles triangles.

**Task 2**. "Is it possible or not to construct any equilateral Diophantine triangle and indicate in its plane such a point, the distance from that to all three vertices of this triangle is expressed in natural numbers?"

We have shown that if does exist such a triangle with sides *a* and such a point that is located respectively on distances *p*, *q* and *k* from its vertices, then occurs the equality:

$$a^4 + p^4 + k^4 + q^4 - a^2p^2 - a^2k^2 - a^2q^2 - p^2k^2 - p^2q^2 - k^2q^2 = 0. \tag{1}$$

A particular solution of (1) is found (it turned out that the found point is located on a circumference circumscribed around this equilateral triangle). At searching other for such a point, it was shown that the interior angle of the triangle, composed by this point and any two vertices of an equilateral triangle, was not found that is commensurate to π.

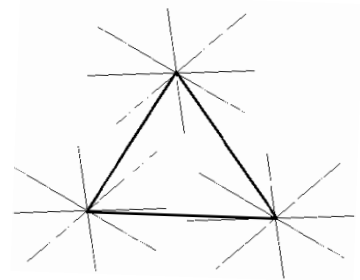

Also is shown that if such a point does exist, then it is impossible to be located on such straight line passing through any vertex of this triangle that makes commensurate to π angle with the sides of this triangle. In addition, we set problematic tasks for the issues mentioned above. A number of tasks on the properties of Diophantine, Bidiophantine, Pseudodiophantine and Pseudobidiophantine geometric figures have been stated and solved.

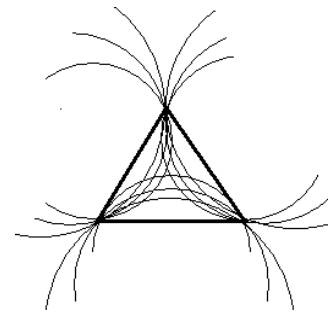

**Basic part**

**Definition 1**. Let's call as Diophantine polygon (polyhedron) to such polygon (polyhedron), the distance between arbitrary teo vertices of that would be expressed by natural number.

**Definition 2**. Let us call as Pseudodiophantine polygon (polyhedron) to such polygon (polyhedron), each side (edge) of that will be expressed by natural number, but length of any diagonal (length of connecting two non-adjacent vertices segment (this segment is not the edge)), is not expressed by rational numbers.

Note that the Steinhaus square is the Pseudodiophantine polygon.

Now let us consider this task of Hugo Steinhaus and assume that does exist such a square and the same point that satisfies the condition of the problem (see Fig. 1). Then we have:

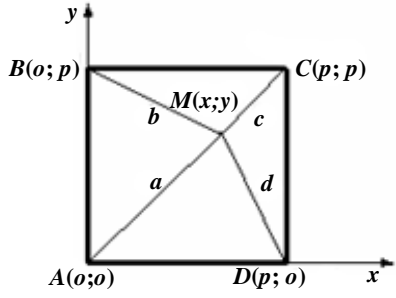

**Fig. 1**

$$\begin{cases} x^2 + y^2 = a^2; \\ x^2 + (p-y)^2 = b^2; \\ (x-p)^2 + (y-p)^2 = c^2; \\ (x-p)^2 + y^2 = d^2. \end{cases} \quad (2)$$

By simple transformations we will obtain:

$$\begin{cases} (p^2 + a^2 - b^2)^2 + (p^2 + a^2 - d^2)^2 = (2pa)^2; \\ (p^2 + b^2 - a^2)^2 + (p^2 + b^2 - c^2)^2 = (2pb)^2; \\ (p^2 + d^2 - a^2)^2 + (p^2 + d^2 - c^2)^2 = (2pd)^2; \\ (p^2 + c^2 - d^2)^2 + (p^2 + c^2 - b^2)^2 = (2pc)^2. \end{cases} \quad (3)$$

Let's assume that from the $\triangle AMB$, $\triangle BMC$, $\triangle CMD$ and $\triangle AMD$ interior angles will not be found angle, the value of that is commensurable to $\pi$.

For this let's consider the several Tasks:

Let us recall the already proved **Lemma 1, Result 1** and **Result 2.**

**Lemma 1.** Let's prove that if $\alpha$ is commensurate to $\pi$ and $\cos\alpha \in Q$, then $\cos\alpha$ would have only and only one one of listed values: $-1; -\frac{1}{2}; 0; \frac{1}{2}; 1$.

Due application of mathematical induction method it is easy to prove that

$$\cos k\alpha = 2^{k-1} \cdot \cos^k \alpha + b_1 \cdot \cos^{k-2} \alpha + \ldots, \quad (4)$$

where $k \in N \setminus \{1\}$; $b_1, b_2, \ldots \in Z \setminus \{0\}$.

From (4) leads that if $\cos\alpha$ is rational, then rational will be also $\cos k\alpha$ ($k \in N$). Now let's say that $\alpha$ is commensurate to $\pi$. I.e. let's say $\alpha = \frac{m}{n}\pi$, where $m \in Z$, $n \in N$ and $\cos\frac{m}{n}\pi = \frac{u}{v}$ ($u, v \in Z \setminus \{0\}$; $(u;v) = 1$). Then in the force of (4)

$$\cos k\frac{m\pi}{n} = 2^{k-1} \cdot \frac{u^k}{v^k} + A_1 \cdot \frac{u^{k-2}}{v^{k-2}} + \ldots, \quad (5)$$

where $k \in N$; $A_1, A_2, \ldots \in Z \setminus \{0\}$.

For $k = n$ by multiplication of both sides on $v^{n-1}$: $A = \frac{2^{n-1} \cdot u^n}{v} + B$ (4), where $A, B \in Z$ and $n, |u|, v \in N$.

If $v$ is different from $2^p$ number, then as it is obvious $\dfrac{2^{n-1} \cdot u^n}{v} \notin Z$, because $(u;v)=1$. This is impossible due (4), thus obligatory $v=2^p$ $p \in Z_0$ (5).

Firstly let's say that $p \in Z_0$ and $p \geq 2$. (6)

In this case let's show that if $\cos n\alpha = \dfrac{a_n}{2^{b_n}}$, where $(a_n;2)=1$ and $n, |a_n|, b_n \in N$, then $b_n > b_{n-1} > \ldots > b_1$ and $(a_1;2)=(a_2;2)=\ldots=(a_n;2)=1$.

By application of mathematical induction method let's prove.

For $n=1$ $\cos \alpha = \dfrac{u}{2^p}$, $(a_1;2)=(u;2)=1$.

For $n=2$ $\cos 2\alpha = \dfrac{u^2 - 2^{2p-1}}{2^{2p-1}}$, $(a_2;2)=(u^2-2^{2p-1};2)=(u^2;2)=1$ and $b_2 = 2p-1 > p = b_1$. Let's assume is valid for $n=k$. I.e.

$$\cos k\alpha = \dfrac{a_k}{2^{b_k}}; \cos(k-1)\alpha = \dfrac{a_{k-1}}{2^{b_{k-1}}}; (a_k;2)=(a_{k-1};2)=1 \text{ and } b_k > b_{k-1} > \ldots > b_1.$$

Let's prove its validity for $n=k+1$

$$\cos(k+1)\alpha = 2\cos k\alpha \cos\alpha - \cos(k-1)\alpha = \dfrac{a_k \cdot u - a_{k-1} \cdot 2^{p+b_k-b_{k-1}-1}}{2^{p+b_k-1}} = \dfrac{a_{k+1}}{2^{b_{k+1}}}.$$

There we have $(a_{k+1};2) = (a_k \cdot u - a_{k-1} \cdot 2^{p+b_k-b_{k-1}-1};2) = (a_k u;2) = 1$ and $b_{k+1} = p + b_k - 1 > b_k$.

We obtain that if $v = 2^p$, where $p \in Z_0 \setminus \{0;1\}$, then $\cos k\alpha$ will not be integer for no one natural number $k$. On the other hand $\cos 2n\alpha = \cos 2n\dfrac{m\pi}{n} = \cos 2\pi m = 1$. I.e. $v = 2^p$ and $p \in \{0;1\}$, from that we obtain the demonstrable.

From this lemma leads that for Diophantine triangles we have:

**Result 1.** If all angles of Diophantine triangle are commensurable to π, then it obligatory will be equilateral.

**Result 2**. If value of one of Diophantine triangle angles is commensurable to π, then it would be only: $\dfrac{\pi}{3}, \dfrac{\pi}{2}$ or $\dfrac{2\pi}{3}$.

**Task 3**. Let's prove that if does exist „square of H. Steinhaus", then from ∠AMB, ∠BMC, ∠CMD and ∠AMD the value of each of them differs from $\dfrac{\pi}{3}$.

Let's assume that the opposite. Without limitation of generality let's say that $\widehat{CMD} = \dfrac{\pi}{3}$. Then accordingly of cosine law from △CMD we will obtain:

$$p^2 = c^2 + d^2 - cd. \tag{7}$$

By introduction in (6) and (7) we will obtain:

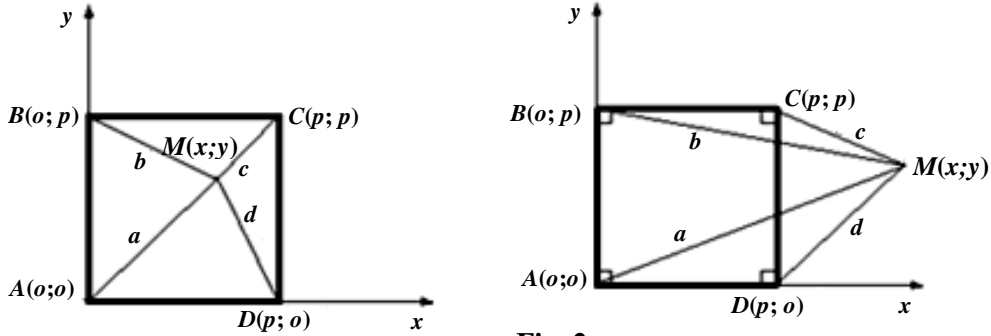

**Fig. 2**

$$\left(c^2 + d^2 - cd + c^2 - d^2\right)^2 + \left(c^2 + d^2 - cd + c^2 - b^2\right)^2 = 4c^2\left(c^2 + d^2 - cd\right) \Leftrightarrow$$

$$\Leftrightarrow \left(2c^2 - cd\right)^2 + \left(2c^2 + d^2 - cd - b^2\right)^2 = 4c^4 + 4c^2d^2 - 4c^3d \Leftrightarrow$$

$$\Leftrightarrow \left(2c^2 + d^2 - cd - b^2\right)^2 = 3(cd)^2.$$

this is impossible, as $a, b, c, d \in N$.

Therefore from $\angle AMB$, $\angle BMC$, $\angle CMD$ and $\angle AMD$ angles value of no one of them is possible to be equal to $\frac{\pi}{3}$.

**Task 4**. Let's prove that if does exist "square of H. Steinhaus", then from $\angle AMB$, $\angle BMC$ $\angle CMD$ and $\angle AMD$ the value of each of them differs from $\frac{2\pi}{3}$.

Let's assume the opposite. Without limitation of generality, let's say that $C\widehat{M}D = \frac{2\pi}{3}$, then from $\angle AMD$, $\angle CMD$ accordingly os cosines law we have:

$$p^2 = c^2 + d^2 + cd. \tag{8}$$

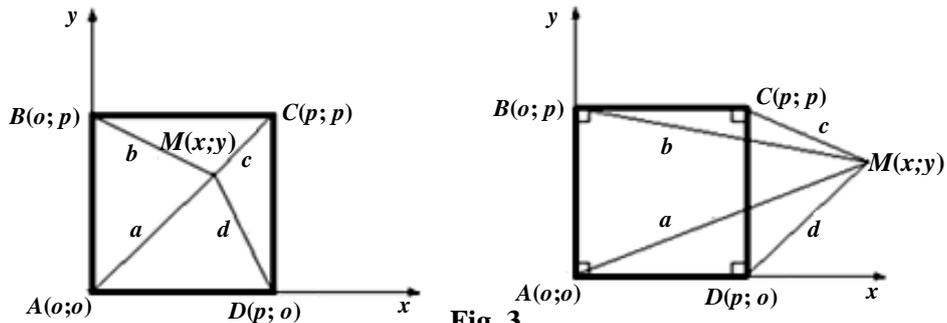

**Fig. 3**

By introduction of (8) in (7) we will obtain:

$$\left(c^2+d^2+cd+c^2-d^2\right)^2+\left(c^2+d^2+cd+c^2-b^2\right)^2=4c^2\left(c^2+d^2+cd\right)\Leftrightarrow$$

$$\Leftrightarrow\left(2c^2+cd\right)^2+\left(2c^2+d^2+cd-b^2\right)^2=4c^4+4c^2d^2+4c^3d\Leftrightarrow$$

$$\Leftrightarrow\left(2c^2+d^2+cd-b^2\right)^2=3(cd)^2, \qquad (9)$$

this is impossible as $a, b, c, d \in N$.

Therefire from angles $\angle AMB$, $\angle BMC$, $\angle CMD$ and $\angle AMD$ the value of no one of them is not possible to be equal to $\dfrac{2\pi}{3}$.

**Task 5.** Lets prove that if does exis "square of H. Steinhaus", then from $\angle AMB$, $\angle BMC$, $\angle CMD$ and $\angle AMD$ the value of each of them differs from $\dfrac{\pi}{2}$.

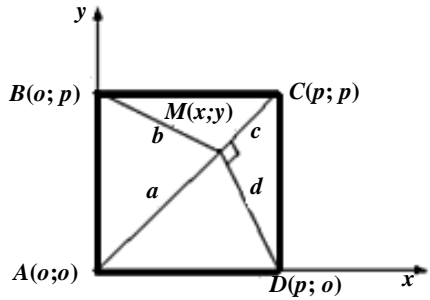 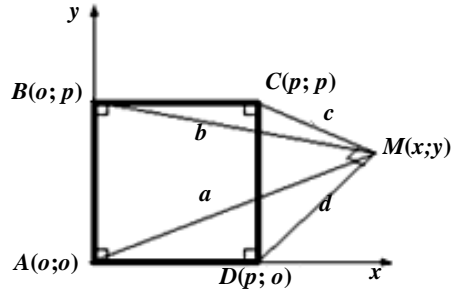

Fig. 4      Fig. 5

There also let's assume the opposite. Without limitation of generality let's say that $\widehat{CMD}=\dfrac{\pi}{2}$. Then accordingly of Pythagorean theorem from $\Delta CMD$ we will obtain:

$$p^2=c^2+d^2. \qquad (10)$$

By introduction of (10) in (6) we will obtain:

$$\left(c^2+d^2+c^2-d^2\right)^2+\left(c^2+d^2+c^2-b^2\right)^2=4c^2\left(c^2+d^2\right)\Leftrightarrow$$

$$\Leftrightarrow\left(2c^2\right)^2+\left(2c^2+d^2-b^2\right)^2=4c^4+4c^2d^2\Leftrightarrow$$

$$\Leftrightarrow 4c^4+\left(2c^2+d^2-b^2\right)^2=4c^4+4c^2d^2\Leftrightarrow$$

$$\Leftrightarrow\begin{bmatrix}2c^2+d^2-b^2=-2cd;\\2c^2+d^2-b^2=2cd.\end{bmatrix}\Leftrightarrow\begin{bmatrix}(c+d)^2+c^2=b^2;\qquad(11)\\(c-d)^2+c^2=b^2.\qquad(12)\end{bmatrix}$$

(11) has not the solution in natural $a, b, c, d$ numbers, as it is reduced to case solved by H. Steinhaus (see Fig. 6, 7).

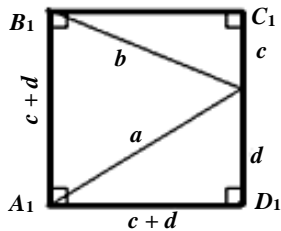 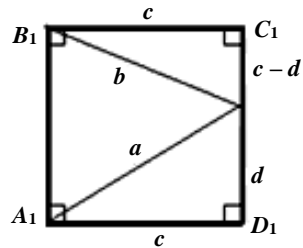

Fig. 6      Fig. 7

(12) also similarly to (11) has not the solution in natural *a, b, c, d* numbers.

I.e. the „square of H. Steinhaus" accordingly of this condition does not exist. Or $A\widehat{M}B \neq \frac{\pi}{2}$, $B\widehat{M}C \neq \frac{\pi}{2}$, $C\widehat{M}D \neq \frac{\pi}{2}$ and $A\widehat{M}D \neq \frac{\pi}{2}$.

From the Task 4, 5 and Task 6 we have:

**Theorem 1**. If does exist. "square of H. Steinhaus", then ∠AMB, ∠BMC, ∠CMD and ∠AMD does not exist angle that is commensurable to π.

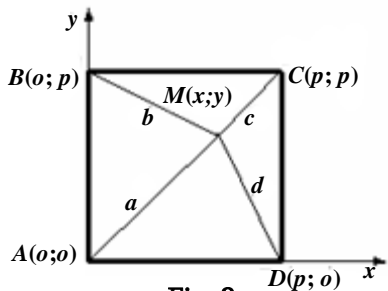

Fig. 8

Let's now show that amongst ∠MAD, ∠MDA, ∠MBC, ∠MBD, ∠MBC, ∠MCB, ∠MBA and ∠MAB does not found the angle value of that also is commensurable to π.

By us in the Task 5 is shown that if α is commensurable to π and $\cos\alpha \in Q$, then $\cos\alpha \in \left\{-1; -\frac{1}{2}; 0; \frac{1}{2}; 1\right\}$. or, in our case $\alpha \in \left\{\pi; \frac{2\pi}{3}; \frac{\pi}{2}; \frac{\pi}{3}; 0\right\}$.

The case of $\alpha \in \left\{0; \frac{\pi}{2}; \pi\right\}$ was considered by H. Steinhaus (see Task 8) and he show that *M* point by such conditions does not exist.

Let's consider the remaining two cases: I. $M\widehat{D}A = \frac{\pi}{3}$, II. $M\widehat{D}A = \frac{2\pi}{3}$.

**I**. If $M\widehat{D}A = \frac{\pi}{3}$, then from ΔAMD $a^2 = p^2 - pd + d^2$ and by its introduction in the third equation of (5) and simplification we will obtain:

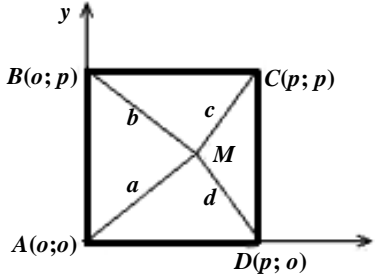

**Fig. 9**

$$(pd)^2 + (p^2 + d^2 - c^2)^2 = (2pd)^2$$

or $(p^2 + d^2 - c^2)^2 = 3(pd)^2$, for that, it is obvious, does not have the solution in natural $p, d, c$ numbers.

**II.** If $M\widehat{D}A = \dfrac{2\pi}{3}$, then from $\triangle MDA$ $a^2 = p^2 + pd + d^2$ by its introduction in the third equation of (5) and simplification we will obtain

$$(pd)^2 + (p^2 + d^2 - c^2)^2 = (2pd)^2$$

or $(p^2 + d^2 - c^2)^2 = 3(pd)^2$ for that, it is obvious, does not have solution in the natural $p, d, c$ numbers.

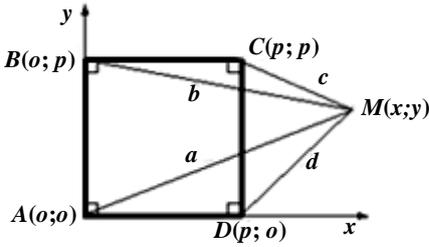

**Fig. 10**

I.e. we obtain that if dies exist such $M$, then it is impossible to be located on such straight lines that passes through any vertex of the square $ABCD$ and makes with any sides of this square the commensurable to $\pi$ angle.

Finally we have the following:

**Theorem 2**. If does exist the "square of H. Steinhaus" ($\square ABCD$), then from the angles $\triangle AMB$; $\triangle BMC$; $\triangle CMD$; $\triangle AMD$; $\triangle MAB$; $\triangle MBA$; $\triangle MBC$; $\triangle MCB$; $\triangle MCD$; $\triangle MDC$; $\triangle MDA$; $\triangle MAD$, would not be found by value commensurable to $\pi$.

Let's now consider the similar to Task ($\square$) task for equilateral triangles Task ($\triangle$).

**Task 7.** Is it possible or not to construct any equilateral Diophantine triangle and indicate in this plane such $M$ point, from that distances to al three vertices of this triangle will be expressed by natural numbers?

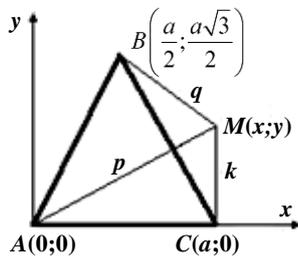

**Fig. 11**

$$\begin{cases} x^2 + y^2 = p^2;\ (x-a)^2 + y^2 = k^2 \\ \left(x - \dfrac{a}{2}\right)^2 + \left(y - \dfrac{a\sqrt{3}}{2}\right)^2 = q^2 \end{cases} \quad (13)$$

From (4.10) we will obtain

$$a^4 - (p^2 + q^2 + k^2)a^2 + (p^4 + q^4 + k^4 - p^2 q^2 - p^2 k^2 - q^2 k^2) = 0 \ . \ (14)$$

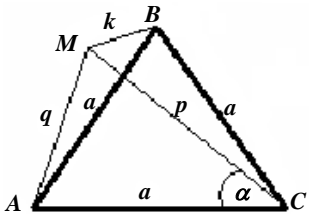

**Fig. 12**

As it was previously mentioned, if $\alpha$ is commensurable to $\pi$ and $\cos\alpha \in Q$, then $\cos\alpha \in \left\{-1; -\dfrac{1}{2}; 0; \dfrac{1}{2}; 1\right\}$, or in our case

$$\alpha \in \left\{\pi; \dfrac{2\pi}{3}; \dfrac{\pi}{2}; \dfrac{\pi}{3}; 0\right\}.$$

Let's say that $M\widehat{C}A = \alpha$. Let us consider each of them:

**I.** If $\alpha = \pi$, then $q = a + p$. By its introduction in (14) and simplification we will obtain:

$$k^4 - 2(a^2 + ap + p^2)k^2 + (a^2 + ap + p^2)^2 = 0 \text{ or}$$
$$k^2 = a^2 + ap + p^2, \qquad (15)$$

all natural solutions of that will be given by:

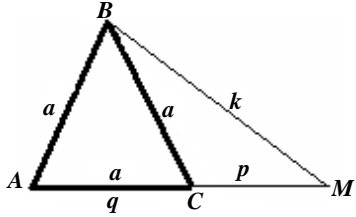

Fig. 13

$$\begin{cases} a = (m^2 - n^2)l \\ p = (m^2 + 2mn)l \\ k = (m^2 + mn + n^2)l \\ q = (2m^2 + 2mn - n^2)l \end{cases} \qquad (16)$$

$$\begin{cases} p = (m^2 - n^2)l \\ a = (m^2 + 2mn)l \\ k = (m^2 + mn + n^2)l \\ q = (2m^2 + 2mn - n^2)l \end{cases} \qquad (17)$$

where $l, m, n \in N$ and $n > m$.

**II.**

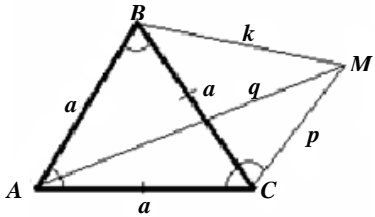

Fig. 14

If $\alpha = \dfrac{2\pi}{3}$, then we have:

$$\begin{cases} k^2 = a^2 - ap + p^2 \\ q^2 = a^2 + ap + p^2. \end{cases} \qquad (18)$$

In the Task 10 we had considered the equation of G.Ch. Poklington, for $k=3$, where by us was shown that (18), that would be reduced to $a^4 + a^2b^2 + p^4 = (kq)^2$, have not the solution in $p, q, k, a$ in natural numbers. I.e. by such conditions such $M$ point does not exist.

**III.** If $\alpha = \dfrac{\pi}{3}$. Then this case is similar to case I.

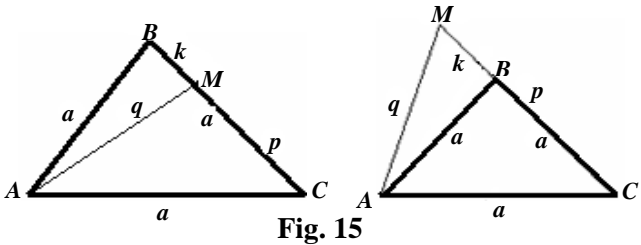

Fig. 15

**IV.** If $\alpha = \dfrac{\pi}{2}$, then in this case we have

$$q^2 = p^2 + a^2. \qquad (19)$$

By introduction of (19) in (14) and simplification we will have

$$k^4 - 2(p^2 + a^2)k^2 + (p^4 - p^2a^2 + a^4) = 0. \quad (20)$$

$$\begin{bmatrix} k^2 = p^2 + a^2 - ap\sqrt{3} & (21) \\ k^2 = p^2 + a^2 + ap\sqrt{3} & (22) \end{bmatrix}$$

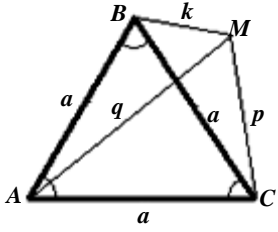

**Fig. 16**

From that, it is obvious that does not have the solution in natural numbers neither (21) as well as (22) and proceeding from this also (20). I.t. *M* point with such conditions does not exist.

**V**. The case $\alpha = 0$ will be reduced to the case I.

I.e. we obtain that desired *M* point is impossible to be located on such straight lines, on that are passing only on one vertex of this equilateral triangle and on containing any side of this triangle straight line represents commensurable to $\pi$ angle. Let's now consider the above mentioned **Task 8**.

**Task 8**. Let's solve kn natural *a, b, c* and *p* numbers the simultaneous equations

$$\begin{cases} a^2 + ab + b^2 = p^2 \\ a^2 - ab + b^2 = c^2. \end{cases} \quad (23)$$

It is obvious that if (23) has the solution in natural *a, b, c* and *p* numbers, then all such solutions would be obtained from such solutions of (23), for that

$$a, b, c, p \in N \text{ and } (a, b) = 1. \quad (24)$$

Indeed, if $(a, b) = k > 1$, then by dividing of each monomials of first and second equations of (23) on $k^2$ we will obtain corresponding to condition (24) condition (23).

Let's solve the (23) by condition $(a, b) = 1$.

Due the multiplication of first and second equations of (23) we will obtain

$$a^4 + a^2b^2 + b^4 = (pc)^2. \quad (25)$$

Firstly let's solve

$$X^2 + XY + Y^2 = Z^2 \quad (26)$$

where

$$X, Y, Z \in N \text{ and } (X; Y; Z) = 1.$$

It is easy to show that all solutions of (26) gives

$$\begin{cases} X = y^2 - x^2 \\ Y = x^2 + 2xy \\ Z = x^2 + xy + y^2 \end{cases}$$

where $x, y \in N, y > x$ and $(x; y) = 1$.

Let's now returen to (25). Let's represent it as

$$(a^2)^2 + a^2b^2 + (b^2)^2 = (pc)^2, \, a, b, c, p \in N \text{ and } (a, b) = 1,$$ at the same time, let's say, *pc* is least. We have

$$\begin{cases} a^2 = y^2 - x^2 \\ b^2 = x^2 + 2xy \\ pc = x^2 + xy + y^2 \end{cases} \Leftrightarrow \begin{cases} a^2 + x^2 = y^2 \\ b^2 = x^2 + 2xy \\ pc = x^2 + xy + y^2 \end{cases} \quad (27)$$

From (27) we have that would ne found such $m, n \in N$;

$$\begin{cases} (m;n) = 1; \\ m > n; \\ mn \equiv 0(\mod 2). \end{cases} \quad (28)$$

For that $\begin{cases} y = m^2 + n^2 \\ \begin{bmatrix} x = m^2 - n^2 \ (*) \\ x = 2mn \ (**) \end{bmatrix} \\ b^2 = x^2 + 2xy \end{cases}$ let's consider two cases:

$$(*)\begin{cases} y = m^2 + n^2; x = m^2 - n^2; \\ b^2 = x^2 + 2xy; \\ b,m,n \in N; (m;n) = 1; m > n; mn \equiv 0(\mod 2). \end{cases} \Rightarrow \begin{cases} b^2 = (m^2 - n^2)^2 + 2(m^4 - n^4) \\ b,m,n \in N; (m;n) = 1; m > n; mn \equiv 0(\mod 2) \end{cases} \Rightarrow$$

$$\Rightarrow \begin{cases} b^2 = (m^2 - n^2)(3m^2 + n^2) = (m^2 - n^2)(4m^2 - (m^2 - n^2)); \\ b,m,n \in N; (m;n) = 1; m > n; mn \equiv 0(\mod 2) \end{cases} \Rightarrow \begin{cases} 3m^2 + n^2 = (2k+1)^2; \\ m^2 - n^2 = (2l+1)^2; \\ m,n,k,l \in N. \end{cases} \Rightarrow$$

$$\Rightarrow \begin{cases} 4m^2 = (2k+1)^2 + (2l+1)^2; \\ m,k,l \in N. \end{cases} \Rightarrow \begin{cases} 4(k^2 + k + l^2 + l) + 2 \equiv 0(\mod 4); \\ k,l \in N. \end{cases}$$

This is impossible.

Let's now consider case (**):

$$\begin{cases} y = m^2 + n^2; \\ x = 2mn; \\ a = m^2 - n^2; b^2 = x^2 + 2xy; \\ b,m,n \in N; (m;n) = 1; m > n; mn \equiv 0(\mod 2) \end{cases} \Rightarrow$$

$$\Rightarrow \begin{cases} b^2 = (2mn)^2 + 2 \cdot 2mn(m^2 + n^2) = 4mn(m^2 + mn + n^2) \\ y = m^2 + n^2; x = 2mn; a = m^2 - n^2; b^2 = x^2 + 2xy; \\ b,m,n \in N; (m;n) = 1; m > n; mn \equiv 0(\mod 2). \end{cases} \quad (29)$$

From (29) with taking into account (28) we have that would be found such $\begin{cases} d, p_1, q \in N; \\ p_1 > q. \end{cases}$ for

that

$$\begin{cases} m = p_1^2 \\ n = q^2 \\ m^2 + mn + n^2 = d^2 \end{cases} \Rightarrow p_1^4 + p_1^2 q^2 + q^4 = d^2. \qquad (30)$$

where $\begin{cases} d, p_1, q \in N; \\ p_1 > q. \end{cases}$

In the case (**) (25) will be as:

$$(p_1^4 - q^4)^4 + (p_1^4 - q^4)^2 (2p_1qd)^2 + (2p_1qd)^4 = (cp)^2 > d^2. \qquad (31)$$

This is impossible, as by our assumption for such type equations $pc$ is least. I.e. also in case (**) (25) have not the solution in natural numbers.

Hence (25) have not the solution in natural numbers.

Let us now show that neither value of $\angle AMB$ and nor $\angle AMC$ is not commensurable to $\pi$. For this it is enough to show for one for them, without the limitation of generality, let's to say that for $\widehat{AMC} = \beta$. There also as in the previous case $\beta \in \left\{\pi; \dfrac{2\pi}{3}; \dfrac{\pi}{2}; \dfrac{\pi}{3}; 0\right\}$.

**I_I.** There also let's consider each of them:

The case $\beta = \pi$ was considered in the case **I**.

**II_I.** If $\beta = \dfrac{2\pi}{3}$, then from $\triangle AMC$ accordingly of cosine law we have:

$$a^2 = p^2 + pq + q^2. \qquad (32)$$

(14) is possible to also written down as:

$$(a^2 - (p^2 + pq + q^2))(a^2 + pq - k^2) = ((p+q)^2 - k^2)(k^2 - (p^2 - pq + q^2)). \qquad (33)$$

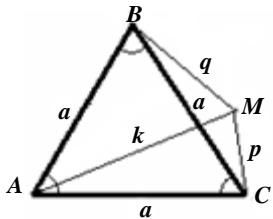

**Fig. 17**

From the Task 8 the system

$$\begin{cases} a^2 - (p^2 + pq + q^2) = 0; \\ k^2 - (p^2 - pq + q^2) = 0. \end{cases}$$

– have not the solution in natural numbers, thus from the (32) and (30) we will obtain that $(p+q)^2 - k^2 = 0$, or $k = p + 1$. This case is considered in the case I.

If $\beta = \dfrac{\pi}{2}$, then from $\triangle AMC$ we have

$$a^2 = p^2 + q^2, \qquad (34)$$

By introduction of that in (34) and simplification we will obtain:

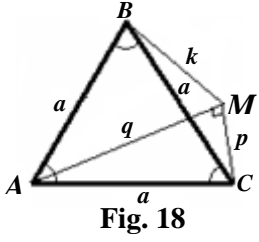

Fig. 18

$$k^4 - 2(p^2 + q^2)k^2 + (p^4 - p^2q^2 + q^4) = 0. \tag{35}$$

From (35) we have

$$\begin{bmatrix} k^2 = p^2 + q^2 - pq\sqrt{3} & (36) \\ k^2 = p^2 + q^2 + pq\sqrt{3} & (37) \end{bmatrix}$$

from that, it is obvious that does not have the solutions in natural numbers neither (36) and nor (37) and hence also nor (35).

I.e. $M$ point by such conditions does not exist.

**IV$_I$.** If $\beta = \dfrac{\pi}{3}$, then from $\triangle AMC$

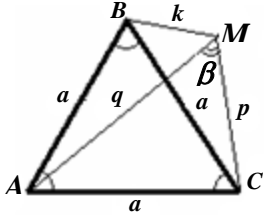

Fig. 19

$$a^2 = p^2 - pq + q^2. \tag{38}$$

(14) also is possible to written down as

$$(k^2 - (p^2 + pq + q^2))(k^2 + pq - a^2) =$$
$$= ((p+q)^2 - a^2)(a^2 - (p^2 - pq + q^2)). \tag{39}$$

From (38) and (39) we have:

$$\begin{cases} a^2 = p^2 - pq + q^2; \\ k^2 = p^2 + pq + q^2. \end{cases} \tag{40}$$

$$\begin{cases} a^2 = p^2 - pq + q^2 \\ k^2 + pq - p^2 + pq - q^2 = 0 \end{cases} \tag{41}$$

(40) has not the solutions in natural numbers accordingly of Task 8.

Let's consider the equation (41).

$$\begin{cases} a^2 = p^2 - pq + q^2 \\ k^2 + pq - p^2 + pq - q^2 = 0 \end{cases} \Leftrightarrow \begin{cases} a^2 = p^2 - pq + q^2 \\ k^2 = (p-q)^2 \end{cases} \Leftrightarrow$$

$$\Leftrightarrow \begin{cases} k + p = q \\ a^2 = k^2 + pk + p^2 \\ a^2 = p^2 - pq + q^2 \end{cases} \Leftrightarrow \begin{cases} q = k + p \\ k^2 + pk + p^2 = p^2 - pq + q^2 \end{cases} \Leftrightarrow$$

$$\Leftrightarrow \begin{cases} q = p + k \\ k^2 + pk + pq - q^2 = 0 \end{cases} \Rightarrow$$

$$\Rightarrow k = \dfrac{-p \pm \sqrt{p^2 - 4pq + 4q^2}}{2} = \dfrac{-p \pm (p - 2q)}{2} \Leftrightarrow$$

$$\Leftrightarrow \begin{bmatrix} k = q - p \\ k = -q \end{bmatrix} \Leftrightarrow q = k + p.$$

$\beta = 0$ is considered in the first case.

Conclusion. In the paper a set of points established by Steinhaus is essentially enlarged and contains Almost everywhere dense set in the plane.